\documentclass[12pt,fleqn, a4]{article}

\usepackage[french]{babel}
\usepackage{inputenc}
\usepackage[T1]{fontenc}

\usepackage{amssymb}
\usepackage{latexsym}
\usepackage{graphics}

\usepackage{graphics}
\usepackage{color}
\newcommand{\colval}{0.3}
\definecolor{colone}{gray}{\colval}

\newcommand{\dcb}{\begin{array}{lll}}
\newcommand{\dce}{\end{array}}
\newcommand{\ebe}{\begin{enumerate}\setlength{\baselineskip}{13pt}\setlength{\parskip}{5pt}}
\newcommand{\dbe}{\end{enumerate}}

\newcommand{\ibegin}{\begin{itemize}\setlength{\baselineskip}{19pt}\setlength{\parskip}{7pt}}
\newcommand{\iend}{\end{itemize}}
\newcommand{\ok}{\rule{4pt}{6pt}}
\newcommand{\desb}{\begin{description}}
\newcommand{\dese}{\end{description}}

\newtheorem{Thm}{Theorem}[section]
\newtheorem {Cor}[Thm]{Corollary}
\newtheorem {definition}[Thm]{Definition}
\newtheorem {pro}{Proposition}[Thm]
\newtheorem {Lemma}[Thm]{Lemma}
\newtheorem {rem}[Thm]{Remark}

\newtheorem {assumption}[Thm]{Assumption}

\newcommand {\bd}{\begin{definition}}
\newcommand {\ed}{\end{definition}}
\newcommand {\bpro}{\begin{pro}}
\newcommand {\epro}{\end{pro}}
\newcommand {\bl}{\begin{Lemma}}
\newcommand {\el}{\end{Lemma}}
\newcommand {\bcor}{\begin{Cor}}
\newcommand {\ecor}{\end{Cor}}
\newcommand {\brem }{\begin{rem} \rm }
\newcommand {\erem }{\end{rem}}
\newcommand{\bethe}{\begin{Thm}}
\newcommand{\ethe}{\end{Thm}}
\newcommand {\bassumption}{\begin{assumption}}
\newcommand {\eassumption}{\end{assumption}}

\def \ind{1\!\!1}
\def\cro#1{\langle #1\rangle}

\begin{document}

\begin{center}
\textbf{\Large An alternative proof of a result of Takaoka}
\end{center}

\begin{center}
Shiqi Song

{\footnotesize Laboratoire Analyse et Probabilités\\
Université d'Evry Val D'Essonne, France\\
shiqi.song@univ-evry.fr}
\end{center}

\

\begin{quote}
\textbf{Abstract.} In Karatzas and Kardaras's paper \cite{KK} it is proved that the NUPBR condition is a property of the local characteristic of the asset process. In Takaoka's paper \cite{takaoka} it is proved that the NUPBR condition is equivalent to the existence of a local martingale deflator. However, the paper \cite{takaoka} founds its proof on Delbaen and Schachermayer's fundamental asset pricing theorem, i.e. the NFLVR condition, which is not a pure property of the local characteristic of the asset process, as shown in \cite{KK}. In this paper we give an alternative proof of the result of \cite{takaoka} which makes use only the properties of the local characteristic of the asset process. This proof complete the program initiated in \cite{KK} which establishes an autonomous theory for the NUPBR condition.  
\end{quote}

\begin{quote}
\textbf{Key words.} NFLVR condition, NUPBR condition, NA condition, local martingale deflator, local characteristic of semimartingale, numéraire portfolio, numéraire change, equivalent $\sigma$-martingale measure, self-financing portfolio

\end{quote}

\begin{quote}
\textbf{JEL classification.} C60, G13
\end{quote}

\begin{quote}
\textbf{MS classification.} 91B70, 60G48
\end{quote}

\

\

\section{Introduction}

We will employ the usual terminology to explain the situation. The notions will be recalled in Section \ref{recall}.

The theory of NFLVR (no free lunch with vanishing risk) as described in the book \cite{DS} is fundamental to the applications of semimartingale models of financial markets. That said, there is little practical tool to check the NFLVR condition and to design equivalent $\sigma$-martingale measures. For this reason, the work of \cite{KK} is of big importance.

Indeed, it is known that the NFLVR condition is the set of two conditions: NA (no arbitrage) and NUPBR (no unbounded profit with bounded risk). It has been shown in \cite{KK} that the NUPBR condition can be checked on the local characteristic of the asset process alone, via the existence of numéraire portfolios. This is a very useful property in practice because most of the models are expressed by their dynamics. Furthermore, the local characteristic is the infinitesimal aspect of a semimartingale, whose study uses only finite dimensional calculus. The proofs of \cite{KK} rely entirely on three ingredients: the semimartingale calculus, the finite-dimensional convex analysis, and measurable selection theorem. In \cite{KK}, it also has been shown that the condition NA can not be determined merely by local characteristic.

Following this work, there was a question:
\begin{quote}
\it Is the NUPBR condition is equivalent to the existence of a local martingale deflator ?
\end{quote}
\cite{K2012} gives an affirmative answer to that question when the dimension $d=1$ and \cite{takaoka} then gives a proof in general case (cf. \cite{Choulli} for a result in continuous case). (In \cite{takaoka} a local martingale deflator is called a strict martingale density. We follow here the terminology of \cite{KK, K2012}). The approaches in \cite{K2012} and in \cite{takaoka} are very different. In this paper we focus on that of \cite{takaoka}.

The idea in \cite{takaoka} is the following. Consider a semimartingale asset processes $S$ which is supposed to be positive and satisfying the NUPBR condition. If the NA condition also holds, we apply the NFLVR theory to get an equivalent local martingale measure, whose density process $\xi$ is then a local martingale deflator for $S$. If the NA condition does not hold, \cite{DS} has proved the existence of a wealth process $\mathsf{V}(H^*)$ which is a maximal element in the set $$
\mathsf{D}=\{f: f\geq 0 \mbox{ and $\exists(H_n)_{\geq 1}$ 1-admissible portfolios such that $H_n^\top\centerdot S\rightarrow f$}\}
$$
Then, making the numéraire change with $\mathsf{V}(H^*)$, the new asset process $\frac{S}{\mathsf{V}(H^*)}$ will satisfy the NFLVR condition with an equivalent local martingale measure. If $\xi$ denotes the density process of that equivalent local martingale measure, the process $\frac{\xi}{W^{H^*}}$ is a local martingale deflator for $S$.

It is to notice nevertheless that \cite{takaoka}'s demonstration does not follow the same schema as in \cite{KK}, because it makes use of the NFLVR condition which depends on factors other than the local characteristic of the asset process. In this note we give a new proof of the result of \cite{takaoka} with an approach in line with \cite{KK}. The proof begins with the following observation : The fundamental theorem of asset pricing of \cite{DS1998} has been proved in two steps (see \cite[Theorem 14.1.1]{DS}). In the first step, one introduces$$
\dcb
\mathsf{K}=\{H^\top\centerdot S: \mbox{ $H$ is 1-admissible for $S$} \}\\
\mathsf{C}=\{f\in\mathsf{L}^\infty(\Omega,\mathcal{A},\mathbb{P}): \exists g\in\mathsf{K}, g\geq f \}\\
\dce
$$
and prove (Kreps-Yan Theorem, see \cite[Theorem 9.1.1 and 14.4.1]{DS}) that there exists an equivalent measure $\mathbb{Q}$ such that 
\begin{equation}\label{KY}
\mathbb{E}_\mathbb{Q}[f]\leq 0 \ \mbox{ for all $f\in\mathsf{C}$.}
\end{equation}
In the second step (\cite[Proposition 14.4.4]{DS}), one applies \cite[The Crucial Lemma 14.3.5]{DS} to prove that, under $\mathbb{Q}$, the opposite of the rate of truncated drift is in the set of barycenters generated by the measures equivalent to the rate of compensation measure. With the close relationship between the equivalent change of the probability measure and the equivalent change of the compensation measure, this property indicates the good way to construct an equivalent probability measure $\mathbb{Q}'$ under which the drift rate of $S$ becomes null. It is to notice that this second step relies uniquely on the local characteristic of $S$, and uses only finite dimensional convex analysis and selection theorem. We see as well that the global property (\ref{KY}) works in correspondence with a similar property at the infinitesimal level (i.e. the assumption of \cite[The Crucial Lemma 14.3.5]{DS}). 

We have then remarked that the property (\ref{KY}) is already satisfied under the initial probability $\mathbb{P}$, if the asset process $S$ is replaced by the numéraire change $\frac{S}{\mathsf{V}(\rho^*)}$, where $\rho^*$ denotes a numéraire portfolio. So, according to the second part of the proof of \cite[Theorem 14.1.1]{DS}, we can assert the existence of an equivalent local martingale measure for $\frac{S}{\mathsf{V}(\rho^*)}$, and therefore a local martingale deflator for $S$ as \cite{takaoka} did.

In fact, we will also give an alternative realization of the idea behind the proof of \cite[Theorem 14.1.1]{DS}.

\section{Semimartingale market model and numéraire portfolio}\label{recall}

\subsection{Setting}

We work on a probability space $(\Omega,\mathcal{A},\mathbb{P})$, where $\mathbb{P}$ denotes a probability defined on the $\sigma$-algebra $\mathcal{A}$ on the space $\Omega$. This probability space is endowed of a filtration $\mathbb{F}=(\mathcal{F}_t)_{t\geq 0}$ of sub-$\sigma$-algebras of $\mathcal{A}$, satisfying the usual conditions. We call an asset process any (multi-dimensional) semimartingale $S=(S_1,\ldots,S_d)$ ($d\geq 1$). For a given asset process $S$ we introduce the family of admissible portfolio $$
\dcb
\mathfrak{A}_S	
=\{H&:& \mbox{$H$ is a $d$-dimensional predictable process,}\\ &&\mbox{integrable with respect to $S$, } \mbox{ and } 1+H^\top\centerdot S> 0\}
\dce
$$
where $H^\top$ denotes the transposition of the vertical vector $H$ and $H^\top\centerdot S$ denotes the stochastic integral of $H$ with respect to $S$. (We refer to \cite{CS1, CS2, jacodbook, JS} for the $d$-dimensional semimartingale integrability. The stochastic integral is denoted by a dot $\centerdot$ and is always supposed null at the origin) We notice that, for any $H\in\mathfrak{A}_S$, there exists a unique self-financing portfolio (cf. subsection \ref{selff}) of the auxiliary asset process $(1,S_1,\ldots,S_d)$, whose corresponding wealth process is given by $1+H^\top\centerdot S$. This quantity will be denoted by $\mathsf{V}(H)$ and called wealth process.

Fix an asset process $S$. Let $T$ to be a $\mathbb{F}$ stopping time representing the time horizon of investments.

\bd
A process $\rho^*\in\mathfrak{A}_S$ will be called a numéraire portfolio for the asset process $S$ on the time horizon $[0,T]$, if, for every $H\in\mathfrak{A}_S$, the quotient process $\frac{\mathsf{V}(H)}{\mathsf{V}(\rho^*)}$ is a supermartingale on $[0,T]$.
\ed

\bd
The asset process $S$ is said to satisfy the NUPBR condition on the time horizon $[0,T]$, if the family of random variables $\{(H^\top\centerdot S)_T: H\in\mathfrak{A}_S\}
$
is bounded in probability.
\ed

\bd
A real valued strictly positive process $Z$ is called a local martingale deflator for the asset process $S$ on the time horizon $[0,T]$, if $Z_T=Z_T\ind_{\{T<\infty\}}+Z_{\infty-}\ind_{\{T=\infty\}}$ exists and strictly positive, and, for any $H\in\mathfrak{A}_S$, $Z(1+H^\top\centerdot S)$ is a local martingale on $[0,T]$.
\ed

\subsection{Strictly positive asset process}\label{s>0}

We consider from now on a $d$-dimensional asset process $S$. We suppose the following assumption in this paper.

\bassumption\label{assumption>0}
The asset process $S=(S_1,\ldots,S_d)$ satisfying $(S_i)_t>0, (S_i)_{t-}>0$, for $1\leq i\leq d$ and $t\geq 0$. 
\eassumption

In the case of this assumption, there exists a $d$-dimensional semimartingale $X=(X_1,\ldots,X_d)$ such that $S_i=(S_i)_0\mathcal{E}(X_i), 1\leq i\leq d$, where $\mathcal{E}$ denotes the stochastic exponential of a semimartingale (see \cite{HWY}). The family $\mathfrak{A}_S$ is, therefore, in one-to-one correspondence with the following class of processes$$
\dcb
\mathfrak{P}_X
=\{\pi&:& \mbox{$\pi$ is a $d$-dimensional predictable process,}\\ &&\mbox{integrable with respect to $X$, } \mbox{ and } \pi^\top\Delta X> -1\}
\dce
$$
where the correspondence is given by $$
\dcb
\pi\in\mathfrak{P}_X &\Longleftrightarrow& \mathcal{E}(\pi\centerdot X)_-\left(\frac{\pi_1}{(S_1)_-},\ldots,\frac{\pi_d}{(S_d)_-}\right) \in\mathfrak{A}_S\\
\\
\frac{1}{\mathsf{V}(H)_-}\left((S_1)_-H_1,\ldots,(S_d)_-H_d\right) \in\mathfrak{P}_X &\Longleftrightarrow& H\in\mathfrak{A}_S \\
\dce
$$
For a corresponding couple $H\in\mathfrak{A}_S$ and $\pi\in\mathfrak{P}_X$, we have $\mathsf{V}(H)=\mathcal{E}(\pi\centerdot X)$. 

\brem
Another implication of Assumption \ref{assumption>0} is the possibility to apply the following result : a positive $\sigma$-martingale is a local martingale (cf. \cite{AS,KK}).  
\erem

\subsection{Local characteristic of the semimartingale $X$}

We refer to \cite{jacodbook, JS} for the notions used in this subsection. Let $\mu_X$ denote the jump measure of $X$ and $\mu_X^p$ denote the compensator (predictable dual projection) of $\mu_X$. Let $X^c$ denote the continuous martingale part of $X$ and $B$ to be the predictable process with finite variation defined in the semimartingale representation of $X$ corresponding to the truncation function $x\ind_{\{|x|\leq 1\}}$ (see \cite[Theorem 2.34]{JS}) :$$
X=X_0+B+X^c+x\ind_{\{|x|\leq 1\}}{_*}(\mu_X-\mu_X^p)+x\ind_{\{|x|> 1\}}{_*}\mu
$$
As in \cite{JS}, we can find four processes $b,c,\nu, G$ satisfying the following conditions : 
\ebe
\item[-]
$b$ is a $d$-dimensional predictable process.
\item[-]
$c$ is a predictable matrix valued process such that, for $t\geq 0$, $c_t$ is a non negative definite symmetric matrix.
\item[-]
$\nu$ is a measure valued predictable process such that, for $t\geq 0$, $\nu_t$ is a Borel measure on $\mathbb{R}^d$ with $\int |x|^2\wedge 1\nu(dx)<\infty$.
\item[-]
$G$ is a real valued increasing predictable process;
\item[-]
The components of $b$ and of $c$ are $G$-integrable. 
\item[-]
It holds :$$
\dcb
\cro{X^c,(X^c)^\top}_t&=&\int_0^t c_sdG_s, \ t\geq 0,\\ 
B_t&=& \int_0^t b_sdG_s, \ t\geq 0,\\
\int_{\mathbb{R}_+\times\mathbb{R}^d} H(s,x)\mu_X^p(ds,dx)&=&\int_0^\infty \left(\int_{\mathbb{R}^d} H(s,x)\nu_s(dx)\right) dG_s,
\dce
$$
where $(X^c)^\top$ denotes the transposition of the vertical vector $X^c$, and $H$ runs over the family of no negative predictable processes.
\dbe
The process $G$ will be called the reference process. When the choice of the process $G$ is fixed, we will call the processes $b$ the truncated drift rate, $c$ the continuous covariation matrix rate, and $\nu$ the compensation measure rate, and we call $(b,c,\nu)$ the local characteristic of $X$ (with respect to $G$).

If $X$ is a special semimartingale, $x\ind_{\{|x|>1\}}{_*}\mu_X^p$ exists. Let $X^v$ denote the drift of $X$ (i.e. $X^v$ is the predictable process with finite variation such that $X-X^v$ is a local martingale). We know that $$
X^v=B+x\ind_{\{|x|>1\}}{_*}\mu_X^p.
$$
We call $x^v=b+\int x\ind_{\{|x|>1\}}\nu(dx)$ the drift rate of $X$.  

Note that the above notation system is valid for any semimartingale.

\subsection{Existence of numéraire portfolios}

We introduce the following set valued processes $\mathfrak{C}^X, \mathfrak{N}^X,\mathfrak{I}^X$ : for $t\geq 0$,
$$
\dcb 
\mathfrak{C}^X_t&=&\{p\in\mathbb{R}^d: \nu_t[p^\top x<-1]=0\},\\
\mathfrak{N}^X_t&=&\{p\in\mathbb{R}^d: p^\top c_t=0, \nu_t[p^\top x\neq 0]=0, p^\top b_t=0\},\\

\mathfrak{I}^X_t&=&\{p\in\mathbb{R}^d: p\notin \mathfrak{N}_t,p^\top c_t=0, \nu_t[p^\top x< 0]=0, p^\top b_t-\int p^\top x\ind_{\{|x|\leq 1\}}\nu_t(dx)\geq 0\}.
\dce
$$
We introduce nextly, for $t\geq 0$, for a $\pi,\rho\in\mathfrak{C}^X_t$, $$
\dcb
\mathfrak{rel}_t(\pi|\rho)=(\pi-\rho)^\top b_t - (\pi-\rho)^\top c_t\rho -\int (\pi-\rho)^\top x\ \frac{1}{1+\rho^\top x}-(\pi-\rho)^\top x\ind_{\{|x|\leq 1\}}\nu_t(dx),\\
\\
\psi_t(\rho)
=\nu_t[\rho^\top x>1]+\left|\rho^\top b_t + \int \rho^\top x(\ind_{\{|x|>1\}}-\ind_{\{|\rho^\top x|>1\}})\nu_t(dx)\right|.
\dce
$$

\bethe(\cite[Theorem 3.15, Theorem 4.12]{KK}, \cite[Lemma 4.1]{K2006})\label{KK}
\ebe
\item
If numéraire portfolios $\rho^*$ exist for $S$ on $[0,T]$, their wealth processes $\mathsf{V}(\rho^*)$ are the same on $[0,T]$. We call this common process on $[0,T]$ the numéraire wealth process and denote it by $\mathsf{V}^*$.
\item
For any $t\geq 0$, $\mathfrak{J}^X_t=\emptyset$ if and only if there exists a unique vector $\rho\in\mathfrak{C}^X_t\cap(\mathfrak{N}^X_t)^\perp$ such that $
\mathfrak{rel}_t(\pi|\rho)\leq 0,\ \forall \pi\in \mathfrak{C}^X_t.
$
We call $\rho$ the pre-numéraire portfolio at $t$.

\item
Suppose that $\mathfrak{J}^X_t=\emptyset$ for $\mathbb{P}\otimes G$ almost all $t\in[0,T]$. We have a predictable process $\rho=(\rho_t)_{t\geq 0}$ such that, for $\mathbb{P}\otimes G$ almost all $t\in[0,T]$, $\rho_t$ is the pre-numéraire portfolio at $t\in[0,T]$. We call $\rho$ a process of pre-numéraire portfolios. 
\item
If the above predictable process $\rho$ of pre-numéraire portfolios satisfies the condition : \linebreak
$\psi(\rho)\centerdot G<\infty$, then $\rho$ is integrable with respect to $X$ and, if we define $$
\rho^*=\mathcal{E}(\rho\centerdot X)_-\left(\frac{\rho_1}{(S_1)_-},\ldots,\frac{\rho_d}{(S_d)_-}\right),
$$ 
$\rho^*$ is a numéraire portfolio for $S$ on $[0,T]$.
\item
Inversely, if a numéraire portfolio $\rho^*$ for $S$ on $[0,T]$ exists, then $\mathfrak{J}^X\neq \emptyset$ has zero $\mathbb{P}\otimes G$ measure on $[0,T]$, and any predictable process $\rho$ of pre-numéraire portfolios satisfies : $\psi(\rho\ind_{[0,T]})\centerdot G<\infty$.
 
\item
A numéraire portfolio $\rho^*$ exists for $S$ on $[0,T]$ such that $\mathsf{V}(\rho^*)_T=\mathsf{V}(\rho^*)_T\ind_{\{T<\infty\}}+\mathsf{V}(\rho^*)_{\infty-}\ind_{\{T=\infty\}}$ exists and finite, if and only if the NUPBR condition holds for $S$ on $[0,T]$.
\dbe
\ethe

\section{Existence of the local martingale deflator}

We suppose in this section that the asset process $S$ satisfies the NUPBR condition on the time horizon $[0,T]$. We assume Assumption \ref{assumption>0}. 

\subsection{Numéraire change}\label{selff}

We need the following known result on numéraire change. Let $V$ be a semimartingale satisfying $V>0$ and $V_->0$. Under this condition, $\frac{1}{V}$ is a semimartingale (see \cite{CKS}). We call the asset process $(\frac{1}{V},\frac{S_1}{V},\ldots,\frac{S_d}{V},1)$ the numéraire change of the asset process $(1,S_1,\ldots,S_d,V)$. A $d+2$-dimensional predictable process $H$ is said a self-financing portfolio for the asset process $(\frac{1}{V},\frac{S_1}{V},\ldots,\frac{S_d}{V},1)$, if $H$ is integrable with respect to the asset process and $$
\left(H^\top(\frac{1}{V},\frac{S_1}{V},\ldots,\frac{S_d}{V},1)^\top\right)_{0+}
=\ \lim_{t\downarrow 0}\left(H^\top(\frac{1}{V},\frac{S_1}{V},\ldots,\frac{S_d}{V},1)^\top\right)_t
$$
exists and
$$
H^\top(\frac{1}{V},\frac{S_1}{V},\ldots,\frac{S_d}{V},1)^\top
=\left(H^\top(\frac{1}{V},\frac{S_1}{V},\ldots,\frac{S_d}{V},1)^\top\right)_{0+}
+H^\top\centerdot(\frac{1}{V},\frac{S_1}{V},\ldots,\frac{S_d}{V},1)^\top.
$$
(The last expression designs a stochastic integral.)

\bl\label{ncvc}
For any bounded self-financing portfolio $H$ with respect to the asset process $(\frac{1}{V},\frac{S_1}{V},\ldots,\frac{S_d}{V},1)$, $H$ is also a self-financing portfolio for the asset process $(1,S_1,\ldots,S_d,V)$. If $$
\lim_{t\downarrow 0}(H^\top\centerdot (1,S_1,\ldots,S_d,V)^\top)_{t}=0,
$$ 
we have$$
H^\top\centerdot (\frac{1}{V},\frac{S_1}{V},\ldots,\frac{S_d}{V},1)^\top
=
\frac{H^\top\centerdot (1,S_1,\ldots,S_d,V)^\top}{V}.
$$
\el

\textbf{Proof.} It is a direct consequence of the integration by parts formula. \ok

By assumption, the NUPBR condition is satisfied by $S$ on $[0,T]$. Let $\rho^*\in\mathfrak{A}_S$ be a numéraire portfolio of $S$ on $[0,T]$ such that $\mathsf{V}(\rho^*)_->0$ on $\mathbb{R}_+$ ($\rho^*$ exists because $\frac{1}{\mathsf{V}^*}$ is a positive supermartingale on $[0,T]$). Consider the process $S'=(\frac{1}{\mathsf{V}(\rho^*)},\frac{S_1}{\mathsf{V}(\rho^*)},\ldots,\frac{S_d}{\mathsf{V}(\rho^*)})$. The components of $S'$ satisfy Assumption \ref{assumption>0}. Let $X'=(X'_0,X'_1,\ldots,X'_d)$ be the semimartingale such that $S'_i=(S'_i)_0\mathcal{E}(X'_i), 0\leq i\leq d$.

\bl\label{piY}
$S'$ is a special semimartingale on $[0,T]$. Moreover, for any $H\in\mathfrak{A}_{S'}$, $1+H^\top \centerdot S'$ is a positive supermartingale on $[0,T]$.
\el

\textbf{Proof.} By definition of numéraire portfolio, the components of $S'$ are positive supermartingales on $[0,T]$, which are special semimartingales.

Let $H=(H_0,H_1,\ldots,H_{d})\in\mathfrak{A}_{S'}$. We suppose firstly that $H$ is locally bounded. We introduce $$
\check{c}_t=
\left\{
\dcb
(H^\top\centerdot S')_{t-} - H_t^\top S'_{t-}, && t>0,\\
\\
0,&& t=0.
\dce
\right.
$$
Then, $(H_0,H_1,\ldots,H_{d},\check{c})$ is a self-financing portfolio for $(\frac{1}{\mathsf{V}(\rho^*)},\frac{S_1}{\mathsf{V}(\rho^*)},\ldots,\frac{S_d}{\mathsf{V}(\rho^*)}, 1)$ such that$$
\lim_{t\downarrow 0}(H_0+\sum_{i=1}^dH_i S_i+\check{c}\mathsf{V}(\rho^*))_t
=\lim_{t\downarrow 0}\frac{1}{\mathsf{V}(\rho^*)_t}(H_0+\sum_{i=1}^dH_i S_i+\check{c}\mathsf{V}(\rho^*))_t=\lim_{t\downarrow 0}(H^\top S'+\check{c})_t=0.
$$ 
According to Lemma \ref{ncvc},$$
\dcb
1+H^\top \centerdot S'
&=&1+(H_0,H_1,\ldots,H_d,\check{c}) \centerdot (\frac{1}{\mathsf{V}(\rho^*)},\frac{S_1}{\mathsf{V}(\rho^*)},\ldots,\frac{S_d}{\mathsf{V}(\rho^*)}, 1)^\top\\

&=&1+\frac{1}{\mathsf{V}(\rho^*)}(H_0,H_1,\ldots,H_d,\check{c}) \centerdot (1,S_1,\ldots,S_d, \mathsf{V}(\rho^*))^\top\\

&=&1+\frac{1}{\mathsf{V}(\rho^*)}(H'^\top\centerdot S+\check{c}\centerdot \mathsf{V}(\rho^*))\\
&&\mbox{ $H'=(H_1,\ldots,H_{d})$}\\

&=&\frac{1}{\mathsf{V}(\rho^*)}(\mathsf{V}(\rho^*)+H'^\top\centerdot S+\check{c}\centerdot \mathsf{V}(\rho^*))\\

&=&\frac{1}{\mathsf{V}(\rho^*)}(1+ (\rho^* +H'+\check{c}\rho^*)^\top\centerdot S)\\

\dce
$$
By definition of $H$, $1+H^\top \centerdot S'>0$. This yields $1+ (\rho^* +H'+\check{c}\rho^*)^\top\centerdot S>0$, i.e. $(\rho^* +H'+\check{c}\rho^*)$ is an element in $\mathfrak{A}_S$. Now, by the definition of numéraire portfolio,$$
1+H^\top \centerdot S' 
=\frac{1}{\mathsf{V}(\rho^*)}(1+ (\rho^* +H'+\check{c}\rho^*)^\top\centerdot S)
=\frac{\mathsf{V}(\rho^* +H'+\check{c}\rho^*)}{\mathsf{V}(\rho^*)}
$$
is a supermartingale on $[0,T]$. 

The above property is proved when $H$ is locally bounded. For a general $H\in\mathfrak{A}_{S'}$, let $\pi$ be its correspondence in $\mathfrak{P}_{X'}$ (see the discussion in subsection \ref{s>0}). For $K>0$, we have $\ind_{\{|\pi|\leq K\}}\pi\in\mathfrak{P}_{X'}$. Let$$
H_K=\ind_{\{|\pi|\leq K\}}\mathcal{E}(\ind_{\{|\pi|\leq K\}}\pi\centerdot X')_-\left(\frac{\pi_1}{(S_1)_-},\ldots,\frac{\pi_d}{(S_d)_-}\right)
$$
Then, $H_K\in\mathfrak{A}_{S'}$, $H_K$ is locally bounded and, therefore, the stochastic exponential\\ 
$$
\mathcal{E}(\ind_{\{|\pi|\leq K\}}\pi\centerdot X')=1+H_K^\top\centerdot S',
$$
and consequently its stochastic logarithm $\ind_{\{|\pi|\leq K\}}\pi\centerdot X'$, is a local supermartingale. Again, 
$$
\ind_{\{|\pi|\leq K\}}\centerdot(1+H^\top S')
=\ind_{\{|\pi|\leq K\}}\centerdot\mathcal{E}(\pi\centerdot X')
=\mathcal{E}(\pi\centerdot X')_-\ind_{\{|\pi|\leq K\}}\pi^\top \centerdot X' 
$$
is a local supermartingale. In other words, $1+H^\top S'$ is a $\sigma$-supermartingale. The lemma is now proved, because a positive $\sigma$-supermartingale bounded at the origin is a supermartingale (cf. \cite[Proposition 11.3]{KK}). \ok

We have immediately the following corollary:

\bcor
$S'$ satisfies the NUPBR condition on $[0,T]$.
\ecor

Without loss of generality, we suppose that $G$ is also a reference process for $S'$.

\subsection{The crucial lemma}

Let $\nu'$ be the rate of compensation measure of $X'$, $b'$ be the rate of truncated drift of $X'$, and $x'^v$ be the drift rate of $X'$. Consider the predictable set $$
D=\{(\omega,t)\in\Omega\times\mathbb{R}_+: \int \ind_{\{|y|>1\}}\nu'_t(dy)(\omega)\neq 0\}
$$

\bl\label{xv=0}
For any predictable process $\pi$ such that $|\pi|<1$, $\ind_{D^c}\pi\in\mathfrak{P}_{X'}$. Consequently, $\ind_{D^c}b'=\ind_{D^c}x'^v=0$ on $[0,T]$.
\el

\textbf{Proof.} For the first part of the lemma, it is enough to note that$$
\nu'_t[\pi^\top x\leq -1]
\leq \nu'_t[|\pi||x|\geq 1]
\leq \nu'_t[|x|> 1]=0, \forall t\notin D.
$$
For the second part, we note that $b'_t=x'^v_t$ for $t\notin D$. Note also that the discussion in subsection \ref{s>0} enable us to interpret Lemma \ref{piY} in term of $X'$. Then, according to Lemma \ref{piY}, by the conclusion of the first part, $\ind_{D^c}\pi^\top x'^v\leq 0$ on $[0,T]$ for all predictable process $\pi$ such that $|\pi|<1$. But this can be possible only when $\ind_{D^c}x'^v=0$ on $[0,T]$. \ok

For $t\in D$ we set $$
\beta_t=\frac{1}{\nu'_t(|x|>1)},\
F_t=(\ind_{\{|x|>1\}}\nu'_t)\star\delta_{\beta b'_t}.
$$
The map $F$ is considered as a measurable map from the set $D$, equipped with the trace $\sigma$-algebra of $\mathcal{P}(\mathbb{F})$ and the trace measure $\mathbb{P}\otimes G$, into the set of finite Borel measures on $\mathbb{R}^d$. 

\bl\label{DScondition}
The map $F$ on the set $D\cap[0,T]$ satisfy the assumption of \cite[The Crucial Lemma 14.3.5]{DS}, i.e., for any predictable process $\pi$ on $D$ such that $F_t(\pi_t^\top z<-1)=0$ for all $t\in D$, we have $\int \pi_t^\top z F_t(dz)\leq 0$ for $t\in D\cap[0,T]$.
\el

\textbf{Proof.} The convolution measure is defined as$$
\int f(z) F_t(dz)=\int f(z) d_z[(\ind_{\{|x|>1\}}\nu'_t)\star\delta_{\beta_tb'_t}]
=\int f(y+\beta_tb'_t) \ind_{\{|y|>1\}} d\nu'_t(dy), \ t\in D.
$$
For any predictable process $\pi$ such that $\pi\in\mathfrak{C}^{X'}$, we have$$
\dcb
\int \pi_t^\top z F_t(dz)
&=&\int \pi_t^\top(y+\beta_tb'_t) \ind_{\{|y|>1\}} d\nu'_t(dy)\\

&=&\pi_t^\top b'_t+\int \pi_t^\top y \ind_{\{|y|>1\}} d\nu'_t(dy)\\
&=&\pi_t^\top x'^v_t, \  t\in D.
\dce
$$
Notice that, for any $K>1$, $\frac{1}{K}\ind_{\{|\pi|\leq K\}}\pi\in\mathfrak{P}_{X'}$. By Lemma \ref{piY}, $\frac{1}{K}\ind_{\{|\pi_t|\leq K\}}\pi_t^\top x'^v_t\leq 0$ for $t\in[0,T]$. Let $K\uparrow\infty$, we prove $\int \pi_t^\top z F_t(dz)\leq 0$ for $ t\in D\cap[0,T]$.

Consider now a predictable process $\pi$ on $D$ such that $F_t(\pi_t^\top z<-1)=0$ for all $t\in D$. Write this condition in the form
$$
\dcb
0&=&F_t(\pi_t^\top z<-1)=\int \ind_{\{\pi_t^\top(y+\beta_tb'_t)<-1\}} \ind_{\{|y|>1\}} d\nu'_t(dy)\\

&=&\int \ind_{\{\pi_t^\top y<-1-\pi_t^\top \beta_tb'_t\}} \ind_{\{|y|>1\}} d\nu'_t(dy)\\

&\geq &\ind_{\{\pi_t^\top \beta_tb'_t>0\}}\int \ind_{\{\pi_t^\top y<-1-\pi_t^\top \beta_tb'_t\}} \ind_{\{|y|>1\}} d\nu'_t(dy)
+\ind_{\{\pi_t^\top \beta_tb'_t\leq 0\}}\int \ind_{\{\pi_t^\top y<-1\}} \ind_{\{|y|>1\}} d\nu'_t(dy)\\

&= &\ind_{\{\pi_t^\top \beta_tb'_t>0\}}\int \ind_{\{\frac{\pi_t^\top y}{1+\pi_t^\top \beta_tb'_t}<-1\}} \ind_{\{|y|>1\}} d\nu'_t(dy)
+\ind_{\{\pi_t^\top \beta_tb'_t\leq 0\}}\int \ind_{\{\pi_t^\top y<-1\}} \ind_{\{|y|>1\}} d\nu'_t(dy)
\dce
$$
Let $$
\dcb
\pi'&=&\ind_D\left(\ind_{\{\pi^\top \beta b'>0\}}\frac{1}{1+\pi^\top \beta b'}+\ind_{\{\pi^\top \beta b'\leq 0\}}\right)\pi,\\
a&=&\left(\ind_{\{\pi^\top \beta b'>0\}}\frac{1}{1+\pi^\top \beta b'}+\ind_{\{\pi^\top \beta b'\leq 0\}}\right),\\
\pi''&=&\frac{a}{2+a|\pi|}\ind_{D}\pi=\frac{\pi'}{2+a|\pi|}.\\

\dce
$$
Then, $$
\dcb
\nu'[\pi''^\top y<-1]
&=&\nu'[\pi''^\top y<-1, |y|\leq 1]+\nu'[\pi''^\top y<-1, |y|> 1]\\
&\leq &\nu'[|\pi''||y|>1, |y|\leq 1]+\nu'[\pi''^\top y<-1, |y|> 1]\\
&\leq &\ind_D\nu'[\frac{a|\pi|}{2+a|\pi|}>1, |y|\leq 1]+\ind_D\nu'[\pi'^\top y<-(2+a|\pi|), |y|> 1]\\
&=&0.
\dce
$$
This nullity shows that $\pi''\in\mathfrak{C}^{X'}$. Consequently, $$
\dcb
0&\geq&\int \pi''^\top_t z F_t(dz)=\frac{a}{2+a|\pi|}\int \pi_t^\top z F_t(dz), \ t\in D\cap[0,T]. \ok
\dce
$$

\subsection{The main result}

We recall the following result from \cite[Théorème II.2]{LM}:

\bl\label{lepingle}
Let $M$ be a local martingale with finite variation. Suppose that the compensator of the increasing process $\sum_{s\leq t}|\Delta_s M|, t\geq 0,$ is bounded. Then, $\mathcal{E}(M)$ has an integrable total variation.
\el

\bethe
Suppose Assumption \ref{assumption>0}. If the NUPBR condition holds for $S$ on $[0,T]$, a local-martingale deflator exists for $S$ on $[0,T]$.
\ethe

\textbf{Proof.} Without loss of generality, we assume that the random variable $G_T$ is bounded. 

Thinks to Lemma \ref{DScondition}, we can apply the \cite[Crucial Lemma 14.3.5]{DS}. We obtain a predictable process $\check{F}$ defined on $D\cap[0,T]$ valued in the set of finite measures on $\mathbb{R}^d$ such that, for $\mathbb{P}\otimes G$ almost all $t\in D\cap[0,T]$,
\ebe
\item[(a')]
$F_t$ is equivalent to $\check{F}_t$, $F_t(\mathbb{R}^d)=\check{F}_t(\mathbb{R}^d)$, and $\|F_t-\check{F}_t\|\leq 1$
\item[(b')]
$\int |y| \check{F}_t(dy)<\infty$ and $\int y \check{F}_t(dy)=0$.
\dbe
Let $\widetilde{F}_t=\check{F}_t\star\delta_{-\beta_t b'_t}$ and $p_t(z)=\frac{d\check{F}_t}{dF_t}(z)$. Then
\ebe
\item[(a)]
For bounded Borel function $f(z)$,  $$
\dcb
\int f(z)\widetilde{F}_t(dz)
&=&\int f(y-\beta_t b'_t)\check{F}_t(dy)
=\int f(y-\beta_t b'_t)p_t(y)F_t(dy)\\
&=&\int f(x+\beta_t b'_t-\beta_t b'_t)p_t(x+\beta_t b'_t)\ind_{\{|x|>1\}}\nu'_t(dx)\\
&=&\int f(x)p_t(x+\beta_t b'_t)\ind_{\{|x|>1\}}\nu'_t(dx).
\dce
$$
This shows that $\widetilde{F}_t$ is equivalent to $\ind_{\{|x|>1\}}\nu'_t$ with density function $p_t(x+\beta_t b'_t)$ which satisfies the inequality :
$$
\dcb
\int |p_t(x+\beta_t b'_t)-1|\ind_{\{|x|>1\}}\nu'_t(dx)=\int |p_t(y)-1|F_t(dy)\leq 1.
\dce
$$
Moreover,$$
\widetilde{F}_t(|x|>1)=\widetilde{F}_t(\mathbb{R}^d)=\check{F}_t(\mathbb{R}^d)=F_t(\mathbb{R}^d)=\nu'_t(|x|>1).
$$
\item[(b)]
$\int |z| \widetilde{F}_t(dz)<\infty$ and $$
\dcb
\int z \widetilde{F}_t(dz)
&=&\int (y-\beta_t b'_t)\check{F}_t(dy)
=\int y \check{F}_t(dy)-\beta_t b'_t \check{F}_t(\mathbb{R}^d)
=-\beta_t b'_t \nu'_t(|x|>1)=-b'_t
\dce
$$
\dbe
Let $U_t(x)$ to be the function $$
U_t(x)=\ind_{\{t\in D\cap[0,T]\}}\ind_{\{|x|>1\}}p_t(x+\beta_t b'_t),\ t\geq 0, x\in\mathbb{R}^d.
$$ 
We have $$
\dcb
(|U-\ind_{D\cap[0,T]}\ind_{\{|x|>1\}}|_*\mu_{X'}^p)_t
&=&\int_0^t\ind_{\{s\in D\cap[0,T]\}}dG_s\int |p_s(x+\beta_s b'_s)-1|\ind_{\{|x|>1\}}\nu'_s(dx)
\leq\int_0^tdG_s\\
\dce
$$
and (cf. \cite{JS} for the notation $\widehat{\cdot}$) $$
\dcb
(U-\ind_{D\cap[0,T]}\ind_{\{|x|>1\}})\widehat{}_s
&=&\ind_{\{s\in D\cap[0,T]\}}\Delta_sG\int (p_s(x+\beta_s b'_s)-1)\ind_{\{|x|>1\}}\nu'_s(dx)\\
&=&\ind_{\{s\in D\cap[0,T]\}}\Delta_sG(\widetilde{F}_s(\mathbb{R}^d)-\nu'(|x|>1))
=0
\dce
$$
Recall that $G_T$ is supposed bounded. So, according to \cite[Chapter II, Theorem 1.33]{JS}, the random function $(U-\ind_{D\cap[0,T]}\ind_{\{|x|>1\}})$ is $(\mu-\mu^p)$ integrable and the stochastic integral $$
M=(U-\ind_{D\cap[0,T]}\ind_{\{|x|>1\}})_*(\mu-\mu^p)
$$ 
is a process with finite variation. Since $\Delta M=(U-\ind_{D\cap[0,T]}\ind_{\{|x|>1\}})$, the previous computations together with Lemmas \ref{lepingle} implies that the exponential martingale $\mathcal{E}(M)$ is positive and uniformly integrable on $[0,T]$ (see also \cite[Lemma 2.3]{K2012}). Applying \cite[Chapter III, Theorem 3.17]{JS}, we can state that, under the new probability $\mathbb{P}'=\mathcal{E}(M)_T\cdot\mathbb{P}$, the compensation measure rate of $X'$ is given by $$
(U-\ind_{D\cap[0,T]}\ind_{\{|x|>1\}}+1)\nu',
$$
and the rate of truncated drift is$$
b'+\int x\ind_{\{|x|\leq 1\}}(U(x)-\ind_{D\cap[0,T]}\ind_{\{|x|>1\}})\nu'(dx)=b'
$$
and the drift rate is$$
\dcb
&&b'+\int x\ind_{\{|x|\leq 1\}}(U(x)-\ind_{D\cap[0,T]}\ind_{\{|x|>1\}})\nu'(dx)+\int x\ind_{\{|x|>1\}}(U(x)-\ind_{D\cap[0,T]}\ind_{\{|x|>1\}}+1)\nu'(dx)\\
&=&b'+\ind_{D\cap[0,T]}\int x\ind_{\{|x|>1\}}(p_t(x+\beta_t b'_t)-1+1)\nu'(dx)
+\ind_{(D\cap[0,T])^c}\int x\ind_{\{|x|>1\}}\nu'(dx)\\

&=&b'\ind_{D\cap[0,T]}+\ind_{D\cap[0,T]}\int x\ind_{\{|x|>1\}}p_t(x+\beta_t b'_t)\nu'(dx)
+b'\ind_{(T,\infty)}+\ind_{(T,\infty)}\int x\ind_{\{|x|>1\}}\nu'(dx)\\
&&\mbox{ according to Lemma \ref{xv=0}}\\

&=&b'\ind_{D\cap[0,T]}+\ind_{D\cap[0,T]}\int z\widetilde{F}(dz)
+b'\ind_{(T,\infty)}+\ind_{(T,\infty)}\int x\ind_{\{|x|>1\}}\nu'(dx)\\

&=&b'\ind_{D\cap[0,T]}-\ind_{D\cap[0,T]}b'
+b'\ind_{(T,\infty)}+\ind_{(T,\infty)}\int x\ind_{\{|x|>1\}}\nu'(dx)\\

&=&b'\ind_{(T,\infty)}+\ind_{(T,\infty)}\int x\ind_{\{|x|>1\}}\nu'(dx)\\\dce
$$
Note that $X'$ is a special semimartingale (cf. Lemma \ref{piY}). The above computation proves that $X'$ is a local martingale on $[0,T]$ under $\mathbb{P}'$. Applying Lemma \ref{piY}, we conclude that, for $H\in\mathfrak{A}_{S'}$, $1+H^\top \centerdot S'$ also is a local martingale on $[0,T]$ under $\mathbb{P}'$. 

Now, let $\xi=\frac{1}{\mathsf{V}(\rho^*)}\mathcal{E}(M)$. We conclude that $\xi$ and $S\xi$ are local martingales under $\mathbb{P}$ on $[0,T]$. For $H\in\mathfrak{A}_S$, $H$ is integrable with respect to $[S,\xi]$ (cf. \cite{CS1}), and consequently, integrable with respect to $$
\xi_-\centerdot S + [S,\xi]=S\xi - S_-\centerdot \xi.
$$
For any $K>0$, we write$$
\ind_{\{|H|\leq K\}}[(1+H^\top\centerdot S)\xi]
=\ind_{\{|H|\leq K\}}(1+H^\top\centerdot S)_-\centerdot \xi+ \ind_{\{|H|\leq K\}}H\centerdot(S\xi - S_-\centerdot \xi)
$$
which is a local martingale under $\mathbb{P}$ on $[0,T]$. This means that $(1+H^\top\centerdot S)\xi$ is a $\sigma$-martingale under $\mathbb{P}$ on $[0,T]$, hence a local martingale, because it is positive. \ok

\

\end{document}